\documentclass{amsproc}

\newtheorem{theorem}{Theorem}[section]
\newtheorem{lemma}[theorem]{Lemma}

\theoremstyle{definition}

\theoremstyle{remark}

\numberwithin{equation}{section}

\begin{document}

\title{Exceptional Vertex Operator Algebras and the Virasoro Algebra}

\author{Michael P. Tuite}
\address{School of Mathematics, Statistics and Applied Mathematics, National University of Ireland, Galway, University Road, Galway, Ireland}
\email{michael.tuite@nuigalway.ie}

\date{November 13, 2008.}

\dedicatory{It is a pleasure to dedicate this paper to Geoffrey Mason on the occasion of his 60th birthday.}


\begin{abstract}
We consider exceptional vertex operator algebras
for which particular Casimir vectors constructed from the primary vectors of lowest conformal weight are Virasoro descendants of the vacuum.  We discuss constraints on these theories that follow from an analysis of appropriate genus zero and genus one two point correlation functions. We find explicit differential equations for the partition function in the cases where the lowest weight primary vectors form a Lie algebra or a Griess algebra.  Examples include the Wess-Zumino-Witten model for Deligne's exceptional Lie algebras and the Moonshine Module. We partially verify the irreducible decomposition of the tensor product of Deligne's exceptional Lie algebras and consider the possibility of similar decompositions for tensor products of the Griess algebra. We briefly discuss some conjectured extremal vertex operator algebras arising in Witten's recent work on three dimensional black holes. 
\end{abstract}

\maketitle

\section{Introduction}
Vertex Operator Algebras (VOAs) are a mathematically rigorous version of chiral conformal field theories and were developed by Borcherds mainly as a tool for understanding Monstrous Moonshine \cite{B}. This work made particular use of the Moonshine module \cite{FLM} whose automorphism group is the Monster sporadic finite simple group. In this paper we consider exceptional VOAs, including the Moonshine Module,  with the property that some quadratic Casimir vectors constructed from the lowest conformal weight primary vectors are Virasoro descendants of the vacuum \cite{Mat,T1}. These conditions give rise to very strong constraints on the VOA which follow from an analysis of appropriate genus zero and genus one two point functions. In the case of VOAs containing a Lie algebra $\frak{g}$ of weight one primaries, then assuming positive rational central charge and integral Kac-Moody level, we find that the genus zero constraints imply that $\frak{g}=A_1,A_2,G_2,D_4,F_4,E_6,E_6,E_8$ known as Deligne's exceptional Lie algebras \cite{D}. The genus one constraints imply that the partition function satisfies an explicit differential equation from which it follows that the VOA is generated by the Lie algebra vectors i.e. the VOA is a level one Wess-Zumino-Witten model for a Deligne exceptional Lie algebra. We also obtain a partial verification of Deligne's original observations for the irreducible decomposition of the tensor product $\frak{g}\otimes \frak{g}$ in this VOA setting. For VOAs for which the lowest primary vectors have weight two and form a Griess algebra, we show that the Virasoro constraints on Casimir vectors imply that the partition function satisfies a specific differential equation and the VOA is generated by the Griess algebra vectors. This restricts a rational central charge to at most nine possible values of which five cases are known to be realizable including the Moonshine module and H\"ohn's Baby Monster VOA \cite{Ho2}. We also discuss possible Deligne-like irreducible tensor product decompositions of the Griess algebra.   We conclude with a brief discussion of similar results for VOAs with primary vectors of lowest weight three, four or five and their relationship to the conjectured extremal VOAs described in Witten's recent work on three dimensional black holes \cite{Wi}. Full details of this work will appear elsewhere \cite{T2}.

\section{Deligne's Exceptional Lie Algebras from Virasoro Constraints}
\subsection{Vertex Operator Algebras.}
We briefly review some notation and properties for Vertex Operator Algebras e.g. \cite{FLM,Ka,MN,LL}. A Vertex Operator Algebra (VOA) consists of a $\mathbb{Z}$-graded vector space $V=\bigoplus_{k\geq 0}V_{k}$ with $V_{0}=\mathbb{C}\mathbf{1}$ for vacuum vector $\mathbf{1}$ and mutually local vertex operators $Y(a,z)=\sum_{n\in \mathbb{Z}}a_{n}z^{-n-1}$ for each $a\in V$
with modes $a_n\in \rm{End}(V)$ such that $a_{-1}\mathbf{1}=a$. There exists a conformal vector $\omega
\in V_{2}$ with $Y(\omega ,z)=\sum_{n\in \mathbb{Z}}L_{n}z^{-n-2}$ 
where $L_{n}$ satisfies the Virasoro algebra of central charge $C$, $L_0$ provides the $\mathbb{Z}$ conformal weight grading and $L_{-1}$ generates translations with $Y(L_{-1}a,z)=\frac{\partial}{\partial z}Y(a,z)$. We consider here VOAs with two further properties (i) 
$V$ is completely reducible with respect to irreducible modules of the  Virasoro algebra and (ii) $V$ is semi-simple \cite{Li}.

Define a symmetric invertible bilinear form $\langle ,\rangle $ on $V$ (with normalization
$\langle \mathbf{1},\mathbf{1}\rangle =1$) where for all $a,b,c\in V$  \cite{FHL},  
\begin{equation}
\langle Y(e^{zL_{1}}(-\frac{1}{z^{2}})^{L_{0}}c,\frac{1}{z})a,b\rangle
=\langle a,Y(c,z)b\rangle.  \label{Yinvariant}
\end{equation}%
For VOAs satisfying (i) and (ii) above, $\langle ,\rangle$ is unique and invertible \cite{Li}. We call such a unique invertible form the Li-Zamolodchikov (Li-Z) metric \cite{MT1} (in deference to ref. \cite{Li} and the equivalent notion in conformal field theory e.g. \cite{BPZ,P}). 
One application of the LiZ metric is in determining the irreducibility of the Verma module $V(C,0)=\bigoplus_{n\geq 0}V_{n}(C,0)$ of Virasoro descendants of the vacuum vector.  For $a,b\in V_{n}(C,0)$ the Gram matrix $M_{n}(C,0)=(\langle a,b\rangle )$ is invertible for all $n$ if and only if $V(C,0)$ is irreducible \cite{KR,Wa}. We recall that there exists a closed formula for the Kac determinant $\det M_{n}(C,0)$ with zeros for central charge \cite{Wa}
\begin{equation}
C_{p,q}=1-6\frac{(p-q)^{2}}{pq},  \label{Cpq}
\end{equation}
where $(p-1)(q-1)=n$ and $(p,q)=1$ for $p,q\geq 2$. In the cases where $C=C_{p,q}$ one constructs the Virasoro minimal model VOA $L(C,0)$ as a quotient of $V(C,0)$ by the radical of $\langle ,\rangle $.   

\subsection{Genus Zero Constraints from Quadratic Casimirs}
Consider a VOA with Li-Z metric for which $d_1=\dim V_1>0$. Defining $[a,b]\equiv a_0b$ for $a,b\in V_1$, then $[a,b]=-[b,a]$ and $[a,b]$ satisfies the Jacobi identity so that $V_1$ defines a Lie algebra with invariant invertible symmetric bilinear form $\langle ,\rangle$ e.g. \cite{Ka,MN}. Furthermore, the modes $\{a_n\}$ of $a\in V_1$ define an affine Kac-Moody algebra. We will denote the Lie algebra by $\frak{g}$ and the corresponding Kac-Moody algebra by $\hat{\frak{g}}$.

Let $\{u_{\alpha }|\alpha =1\ldots d_1\}$ and  
$\{u^{\beta }|\beta =1\ldots d_1\}$ denote a $\frak{g}$ basis and LiZ dual basis, respectively i.e. $\langle u^{\alpha },u_{\beta}\rangle =\delta _{\beta }^{\alpha }$. We define the quadratic Casimir vectors 
\begin{equation}
\lambda^{(n)}=u^{\alpha}_{1-n}u_{\alpha}\in V_n,\label{eq:lambdan}
\end{equation}
where $\alpha$ is summed. It immediately follows that $\lambda^{(0)}=-d_1\mathbf{1}$ and 
$\lambda^{(1)}=0$. Furthermore
\begin{equation}
L_m\lambda^{(n)}=(n-1)\lambda^{(n-m)},\quad m>0\label{eq:LMlambda}
\end{equation}
implying the following lemma concerning Virasoro descendants of the vacuum \cite{Mat}:
\begin{lemma}\label{lemma:lambda}
If $\lambda^{(n)}\in V_n(C,0)$ and $\det M_{n}(C,0)\neq 0$ then $\lambda^{(n)}$ is uniquely determined.
\end {lemma}
Thus if $\lambda^{(2)}\in V_2(C,0)$ then  $\lambda^{(2)}=\kappa\omega$ for some $\kappa\in \mathbb{C}$ so that $\langle \omega ,\lambda^{(2)}\rangle =\kappa\frac{C}{2}$ for $C\neq C_{2,3}=0$ where $\det M_2(C,0)=0$. But (\ref{Yinvariant}) and (\ref{eq:LMlambda}) imply
$\langle \omega,\lambda ^{(2)}\rangle =
\langle \mathbf{1},L_{2}\lambda ^{(2)}\rangle =-d_1$ so that 
\begin{equation}
\lambda^{(2)}=-\frac{2d_1}{C}\omega.\label{eq:lambda2}
\end{equation}
For $\frak{g}$ simple, this is the standard Sugawara construction for $\omega$ e.g. \cite{Ka}. Similarly, if  $\lambda^{(4)}\in V^{(4)}(C,0)$ for $C\neq C_{2,3}=0,C_{2,5}=-22/5$ (for which $\det M_4(C,0)=0$) then \cite{MMS,T1}
\begin{equation}
\lambda^{(4)}=-\frac{3d_1}{C\left(5C+22\right)}\left[4L^2_{-2}\mathbf{1}%
+\left( 2+C\right) L_{-4}\mathbf{1}\right].  \label{eq:lambda4}
\end{equation}
We also note that (\ref{eq:lambda4}) implies (\ref{eq:lambda2}) using (\ref{eq:LMlambda}).

We next consider the constraints on $\frak{g}$ if either (\ref{eq:lambda2}) or (\ref{eq:lambda4}) hold. In order to do this it is useful to introduce the correlation function  
for $a,b \in \frak{g}$
\begin{eqnarray}
F(a,b;x,y)=\langle a,Y(u^{\alpha},x)Y(u_{\alpha},y)b\rangle,\label{eq:Fab}
\end{eqnarray}
with $\alpha$ summed. Using the associativity property of vertex operators \cite{FLM,Ka,MN} and (\ref{eq:lambdan}) we find
\begin{equation}
F(a,b;x,y)=\sum_{n\ge 0} \langle a,o(\lambda^{(n)})b\rangle (x-y)^{n-2},\label{eq:Fablambda}
\end{equation}  
for \lq zero mode\rq\  $o(a)\equiv a_{n-1}:V_k\rightarrow V_k$ for $a\in V_n$. $F(a,b;x,y)$ may be formally expanded in $x,y$ in a number of ways so as to obtain a rational expression in terms of the LiZ metric 
$\langle a,b\rangle$ and the Lie algebra Killing form 
$K(a,b)\equiv Tr_{\frak{g}}(a_0b_0)$ as follows \cite{T1}:
\begin{theorem}\label{theorem:FabG}
$F(a,b;x,y)$ is given by the rational function
\begin{equation*}
F(a,b;x,y)=\frac{G(a,b;x,y)}{x^{2}y^{2}(x-y)^{2}},
\end{equation*}
where
\begin{eqnarray*}
G(a,b;x,y)=&&-\left[d_1x^{2}y^{2}+2xy(x-y)^{2}+(x-y)^{4}\right]\langle a,b\rangle\\
 &&+xy(x-y)^{2}K(a,b).
\end{eqnarray*}
\end{theorem}
Since $F(a,b;x,y)$ is rational, $x,y$ can be considered as points on the genus zero Riemann sphere so that $F(a,b;x,y)$ is referred to as a genus zero correlation function. Later on we will construct correlation  functions associated with a genus one torus. 

Suppose that $\lambda^{(2)}$ is a vacuum Virasoro descendant so that (\ref{eq:lambda2}) holds. Then (\ref{eq:Fablambda}) results in a further constraint on $F(a,b;x,y)$ which determines the Killing form: \cite{T1}
\begin{theorem}\label{theorem:Killing}
If $\lambda^{(2)}\in V_2(C,0)$ and $C\ne 0$ then 
\begin{equation}
K(a,b)=-2\langle a,b\rangle \left( \frac{d_1}{C}-1\right).\label{eq:KillingForm}
\end{equation}
\end{theorem}
Since the Li-Z metric is non-degenerate, it immediately follows from
Cartan's criterion that $\frak{g}$ is semi-simple for $d_1\neq C$ and is solvable for $d_1=C$. For $d_1\neq C$ this implies the Kac-Moody algebra decomposes as 
\begin{equation}
\hat{\frak{g}}=\hat{\frak{g}}_1^{(k_1)}\oplus \hat{\frak{g}}_2^{(k_2)}\oplus \ldots \oplus \hat{\frak{g}}_r^{(k_r)},  \label{V1simple}
\end{equation}
where $\hat{\frak{g}}_i^{(k_i)}$ is a simple Kac-Moody algebra of level 
$k_{i}=-\frac{1}{2}\langle \mathbf{\alpha }_{i},\mathbf{\alpha }_{i}\rangle$ for a long root $\mathbf{\alpha }_{i}$ and with dual Coxeter number 
\begin{equation}
h_{i}^{\vee }=k_{i}(\frac{d_1}{C}-1).  \label{hoverk}
\end{equation}

Suppose further that $\lambda^{(4)}$ is a vacuum Virasoro descendant so that (\ref{eq:lambda4}) holds. This results in a further constraint from (\ref{eq:Fablambda}) on $F(a,b;x,y)$ so that \cite{T1}
\begin{theorem}\label{theorem:d(C)}
If  $\lambda^{(4)}\in V_{4}(C,0)$ and $C\ne 0,-\frac{22}{5}$ then  
\begin{equation}
d_1=\frac{C(5C+22)}{10-C}.\label{eq:d(C)}
\end{equation}
\end{theorem}  
For integral $d_1>0$ there are 21 positive rational values of $C$ given by 
\begin{equation}
C={\frac {2}{5}},1,2,{\frac {14}{5}},4,5,{\frac {26}{5}},6,{\frac {32}{5}},{\frac 
{34}{5}},7,{\frac {38}{5}},8,{\frac {41}{5}},{\frac {42}{5}},{\frac {
44}{5}},9,{\frac {46}{5}},{\frac {47}{5}},{\frac {48}{5}},{\frac {49}{
5}}.\label{eq:gCvalues}
\end{equation}
There are also 21 negative rational solutions with
\begin{equation}
C^{\ast }=-2d_1/C,\label{eq:C*} 
\end{equation}
where $C\leftrightarrow C^\ast$ interchanges the roots of the quadratic (\ref{eq:d(C)}). Applying Theorem \ref{theorem:Killing}, Cartan's condition implies that $\frak{g}$ must be simple of some Kac-Moody level $k_1$. Restricting $k_1$ to be integral\footnote{Note that this condition was mistakenly omitted in ref. \cite{T1}} (for example, if $V$ is assumed to be $C_2$-cofinite \cite{DM}) we find 
\begin{theorem}\label{theorem:Deligne}
Suppose $\lambda^{(4)}\in V_{4}(C,0)$, $C\ne 0,-\frac{22}{5}$ positive rational and that the simple component Kac-Moody levels are integral. Then $\frak{g}$ is one of the simple Lie algebras  
$A_1,A_2,G_2,D_4,F_4,E_6,E_6,E_8$ for central charge $C=1,2,\frac{14}{5},4,\frac{26}{5},6,7,8$ respectively with level one Kac-Moody algebra $\hat{\frak{g}}^{(1)}$  and dual Coxeter number $h^{\vee }(C)=\frac{d_1}{C}-1=\frac{12+6C}{10-C}$. 
\end{theorem}

A similar result is obtained in \cite{MMS} based on a number of stronger assumptions. The simple Lie algebras appearing in Theorem \ref{theorem:Deligne} are known as Deligne's exceptional Lie algebras \cite{D}.  Deligne's list has also been noted in some specific lattice VOA constructions \cite{Hu}.

\subsection{Deligne's Exceptional Lie Algebras.}
Deligne observed that many properties of the Lie algebras $\frak{g}=A_1,A_2,G_2,D_4,F_4,E_6,E_6,E_8$ can be described universally in terms of a single parameter \cite{D}. Choosing the parameter to be 
$h^{\vee }$ \cite{C} then $\dim\frak{g}$ is observed to be 
\begin{equation}
\dim\frak{g}=\frac{2(5h^{\vee }-6)(h^{\vee }+1)}{h^{\vee }+6},\label{eq:Vogel}
\end{equation}%
known as the Vogel formula. Thus $\dim\frak{g}$ precisely agrees with $d_1$ of (\ref{eq:d(C)}) for $h^{\vee }(C)$ of Theorem \ref{theorem:Deligne} so that the central charge $C$ provides an alternative parameterization for Deligne's observations.
  
Deligne, and later others, also found universal rational formulas for the dimension of the irreducible components of the tensor product $\bigotimes\limits^{n}\frak{g}$ for $n=2,3,4$ \cite{D,DdeM,CdeM,C} and for parts of arbitrary tensor products \cite{LM}. For example, the symmetric part of $\frak{g}\otimes \frak{g}$ universally decomposes as 
\begin{equation}
\rm{Sym}\left(\frak{g}\otimes \frak{g}\right) =1\oplus Y\oplus Y^{\ast },\label{eq:Symg2}
\end{equation}
where $Y$ and $Y^{\ast}$ denote two irreducible representations of $\frak{g}$ of dimension (expressed here in terms of $C$)  
\begin{equation}
\dim Y=\frac{5(5C+22)(C-1)(C+2)^2}{2(C-22)(C-10)}.\label{eq:dimY}
\end{equation}
Since $\dim{\rm Sym}\left(\frak{g}\otimes \frak{g}\right)=\frac{1}{2}d_1(d_1+1)$
and $d_1(C)=d_1(C^\ast)$ (from (\ref{eq:C*})) it follows that $\dim Y^{\ast }=\dim Y(C^{\ast })$. 

\subsection{Genus One Constraints from Quadratic Casimirs}
We next consider the constraints on the genus one partition function $Z(q)$ that follow from the Virasoro descendant condition (\ref{eq:lambda4}). We will show that in this case, $Z(q)$ satisfies a second order differential equation and hence is uniquely determined. As a consequence, we prove that $V$ is a level one Wess-Zumino-Witten VOA for $V_1=\frak{g}$ a Deligne exceptional series. We will also partly explain the irreducible decomposition of $\rm{Sym}\left( \frak{g}\otimes \frak{g}\right)$ in (\ref{eq:Symg2}). Full details will appear elsewhere \cite{T2}.
   
Define the standard genus one partition function of a VOA $V$ by the trace function
\begin{equation}
Z(q)=Tr(q^{L_{0}-C/24})=q^{-C/24}\sum_{n\geq 0}\dim V_n q^n,\label{eq:Z(q)}
\end{equation}
where $q$ is, at this stage, a formal parameter. In order to analyze such functions and to develop  a theory of genus one n-point correlation functions, Zhu  \cite{Z} introduced an isomorphic VOA with \lq square bracket\rq\  vertex operators   
\begin{equation}
Y[a,z]\equiv Y(e^{zL_{0}}a,e^{z}-1)=\sum\limits_{n\in \mathbb{Z}}a_{[n]}z^{-n-1},\label{eq:square}
\end{equation}%
and Virasoro vector $\tilde{\omega}=\omega -\frac{C}{24}\mathbf{1}$  with modes $\{L[n]\}$. In particular, $L[0]$ defines an alternative $\mathbb{Z}$ grading on $V$ with $V=\bigoplus_{k\geq 0}V_{[k]}$.

We define the genus one 1-point correlation function for $a\in V$ by 
\begin{equation}
Z(a,q)=Tr(o(a)q^{L_{0}-C/24}),\label{eq:1pt}
\end{equation}
for zero mode $o(a)=a_{n-1}$ for $a\in V_n$. Thus, $Z(\tilde{\omega},q)=Tr((L_{0}-\frac{C}{24}%
)q^{L_{0}-C/24})=q\frac{\partial }{\partial q}Z(q)$.
We define the 2-point correlation function in terms of 1-point functions by
\begin{eqnarray}
Z((a,x),(b,y),q) &=&Z(Y[a,x]Y[b,y]\mathbf{1},q)\label{eq:2pt1} \\
&=&Z(Y[a,x-y]b,q),\label{eq:2pt2}
\end{eqnarray}%
where the second identity follows from associativity \cite{MT2}.

We now consider the constraints on the partition function $Z(q)$ arising from the Virasoro descendant condition (\ref{eq:lambda4}). We define square bracket quadratic Casimir vectors $\lambda^{[n]}=u^{\alpha}_{[1-n]}u_{\alpha}\in V_{[n]}$. Since the square bracket VOA with Virasoro vector $\tilde\omega$ is isomorphic to the original VOA with Virasoro vector $\omega$ it follows that $\lambda^{(n)}\in V_{n}(C,0)$ if and only if $\lambda^{[n]}\in V_{[n]}(C,0)$.  We define a genus one analogue of (\ref{eq:Fab}) given by the 2-point function $Z((u^{\alpha },x),(u_{\alpha },y),q)$. Furthermore, associativity (\ref{eq:2pt2}) implies the genus one analogue of (\ref{eq:Fablambda})
\begin{eqnarray}
Z((u^{\alpha },x),(u_{\alpha },y),q)=\sum_{n\geq 0}Z(\lambda^{[n]},q)(x-y)^{n-2}.\label{eq:Zlambda}
\end{eqnarray}
Zhu proved a reduction formula (Proposition 4.3.2 \cite{Z}) where an $n$-point function is expanded in terms of $n-1$ point functions with coefficients given by elliptic Weierstrass functions $P_{n+2}(z,q)=%
\frac{(-1)^{n}}{(n+1)!}\frac{\partial^{n}}{\partial z^{n}}P_2(z,q)$ for
\begin{equation}
P_2(z,q)=\frac{1}{z^{2}}+\sum_{k\geq 2}E_{k}(q)(k-1)z^{k-2},\label{eq:P2}
\end{equation}
with modular Eisenstein series $E_{k}(q)=0$ for $k$ odd and for $k\ge 2$ even
\begin{equation}
E_{k}(q)=-\frac{B_{k}}{k!}+\frac{2}{(k-1)!}
\sum_{n\geq 1}\frac{n^{k-1}q^{n}}{1-q^{n}}, \label{eq:En}
\end{equation}
and $B_{k}$ is the $k$th Bernoulli number. The parameter $z$ can thus be considered as a point on the genus one torus with modular parameter $\tau$ where $q=\exp(2\pi i \tau)$ for $|q|<1$. Applying Zhu reduction to the LHS of (\ref{eq:Zlambda}) leads to
\begin{eqnarray*}
Z((u^{\alpha },x),(u_{\alpha },y),q)= Tr(o(u^{\alpha })o(u_{\alpha })q^{L_{0}-C/24})
-d_1P_2(x-y,q)Z(q).
\end{eqnarray*}
Comparing the  $(x-y)^{2}$ term in this expression and the RHS of (\ref{eq:Zlambda}) we find
\begin{equation}
Z(\lambda^{[4]},q)=-3d_1E_{4}(q)Z(q).\label{eq:Zlambda4}
\end{equation} 
If $\lambda^{[4]}$ is a square bracket Virasoro descendant then $Z(\lambda^{[4]},q)$ can be evaluated using an alternative Virasoro Zhu reduction. Then we find:
\begin{theorem}\label{theorem:Zdeqn}
If $\lambda^{(4)}\in V_{4}(C,0)$ for $C\ne 0,-\frac{22}{5}$ then $Z(q)$ is the unique solution with leading form $Z(q)=q^{-C/24}(1+O(q))$ to the differential equation
\begin{equation}                                                        
\left[\left(q\frac{\partial}{\partial q} \right)^2+2E_2(q)q\frac{\partial}{\partial q}
-\frac{5}{4}C(C+4)E_{4}(q)\right]Z(q)=0.\label{eq:Zdeqn} 
\end{equation}
$Z(q)$ is convergent for $0<|q|<1$.
\end{theorem}

\textbf{Proof.} 
$\lambda^{(4)}\in V_{4}(C,0)$ if and only if $\lambda^{[4]}\in V_{[4]}(C,0)$. Thus   
$\lambda^{[4]}=\frac{3d_1}{C(5C+22)}(4L_{[-2]}^{2}\mathbf{1}+(2+C)
L_{[-4]}\mathbf{1})$.  Zhu reduction implies $Z(L_{[-4]}\mathbf{1},q)=0$ and
\begin{eqnarray*}
Z(L_{[-2]}^{2}\mathbf{1},q)=\left[\left(q\frac{\partial}{\partial q} \right)^2+2E_2(q)q\frac{\partial}{\partial q}+\frac{C}{2}E_{4}(q)\right] Z(q),\\
\end{eqnarray*}
from which (\ref{eq:Zdeqn}) follows on using (\ref{eq:Zlambda4}). (\ref{eq:Zdeqn}) has a regular singular point at $q=0$ with indicial roots $-C/24$ and $(C+4)/24$. Thus there exists a unique solution with leading form $Z(q)=q^{-C/24}(1+O(q))$ where $Z(q)$ converges for $0<|q|<1$.   \qed

Much as in the approach taken by Zhu for $C_2$-cofinite theories \cite{Z}, we note that the space of solutions to the differential equation (\ref{eq:Zdeqn}) is modular invariant as is discussed further in \cite{Mas}. 
 
An immediate consequence of Theorem \ref{theorem:Zdeqn} is:
\begin{theorem}\label{theorem:ZWZW}
$V$ is generated by $\frak{g}$ i.e. $V$ is a Wess-Zumino-Witten VOA. 
\end{theorem}

\textbf{Proof.} Consider the subVOA $V'\subseteq V$ generated by $\frak{g}$ with partition function $Z'(q)$. Clearly $\lambda^{(2)},\lambda^{(4)}\in V'$. If $\lambda^{(4)}$ is a vacuum descendant of $V$ then $\lambda^{(2)}=-\frac{2d_1}{C}\omega$ and hence $\omega\in V'$. Thus $\lambda^{(4)}$ is a vacuum descendant of $V'$ also and  it follows that $Z'(q)$ obeys (\ref{eq:Zdeqn}). Hence $Z'(q)=Z(q)$ and so $V'=V$.\qed

Combining Theorems \ref{theorem:Deligne} and \ref{theorem:ZWZW} we thus find:
\begin{theorem}
Suppose that $\lambda^{(4)}\in V_{4}(C,0)$, $C\ne 0,-\frac{22}{5}$ is positive rational and the simple component Kac-Moody levels are integral.
Then $V$ is the Kac Moody level 1 WZW VOA for $\frak{g}=A_{1},A_{2},G_{2},D_{4},F_{4},E_{6},E_{7},E_{8}$.
\end{theorem}

\subsection{$\frak{g}\otimes\frak{g}$ Irreducible Structure.}
Consider the vector space of weight two vectors $V_2$. We may decompose this according to its Virasoro Verma module structure as  
\begin{equation}
V_2=\mathbb{C}\omega\oplus L_{-1}\frak{g}\oplus P_2,\label{eq:V_2}
\end{equation}
where $P_2$ is the space of weight two primary vectors i.e. $L_n a=0$ for all $n>0$ for $a\in P_2$. Thus $\dim V_2 =1+d_1+p_2 $ with $p_2=\dim P_2$. Substituting $Z(q)=q^{-C/24}\sum_{n}\dim V_n q^n$ into the differential equation (\ref{eq:Zdeqn}) we may solve recursively for $\dim V_n$ as a rational function in $C$. (Since $\dim V_n$ is positive integral, this restricts the possible rational values of $C$ in (\ref{eq:gCvalues}) even before we invoke Cartan's criterion following Theorem \ref{theorem:Killing}). In particular, for $d_1=\dim \frak{g}$ we recover (\ref{eq:d(C)}) whereas solving for $\dim V_2$ results in
\begin{eqnarray*}
p_2 =\frac{5( 5C+22)(C-1)(C+2)^2}{2(C-22)(C-10)}.
\end{eqnarray*}
This is precisely Deligne's formula for the dimension of the irreducible
representation $Y$ in (\ref{eq:dimY}). 

From Theorem \ref{theorem:ZWZW} we know $V$ is generated by $\frak{g}$ so that
\begin{equation*}
V_2=\mathbb{C}[u_{-2}^{\alpha }\mathbf{1},u_{-1}^{\alpha }u_{-1}^{\beta }\mathbf{1}].
\end{equation*}
Clearly the Virasoro descendants $L_{-1}\frak{g}$ and the antisymmetric part of $\mathbb{C}[u_{-1}^{\alpha }u_{-1}^{\beta }\mathbf{1}]$ (using the Kac-Moody Lie algebra) lie in $\mathbb{C}[u_{-2}^{\alpha }\mathbf{1}]$ whereas the Virasoro vector $\omega$ and the primary vectors $P_2$ lie in the symmetric part of $\mathbb{C}[u_{-1}^{\alpha }u_{-1}^{\beta }\mathbf{1}]$. This concurs with the occurrence of the representation $Y$ of dimension $p_2$ in (\ref{eq:Symg2}).

We briefly sketch how to prove that $P_2$ is an irreducible representation of $\frak{g}$ by use of character theory \cite{T2}. The main idea is to use Zhu reduction formulas for orbifold trace functions which include a Lie group element $g=\exp(a_0)$ generated by $a_0$ for $a\in \frak{g}$ \cite{MTZ}. Thus we consider genus one trace functions of the form 
\begin{eqnarray*}
	Z(g,q)&=&Tr(gq^{L_0-C/24})\\
	&=&q^{L_0-C/24}\left[1+\chi_1(g)q+(1+\chi_1(g)+\chi _{2}(g))q^2\ldots  \right],
\end{eqnarray*}
$\chi_{1}(g)$ and $\chi _{2}(g)$  are the characters for $\frak{g}$ and  $P_2$ respectively following the decomposition of (\ref{eq:V_2}). 
By considering an appropriate 2-point function, a generalized version of Theorem \ref{theorem:Zdeqn} can be found leading to a differential equation for $Z(g,q)$ involving \lq twisted\rq\ Eisenstein series \cite{MTZ}. Analyzing this equation results in
  
\begin{theorem}
The characters $\chi_{1}(g),\chi _{2}(g)$ are related
as follows 
\begin{eqnarray*}
\chi _{2}(g) &=&\frac{1}{C-22}\big{[}5C+22+3(2+C)\chi _{1}(g) \\
&&+\frac{1}{2}(C-10)( \chi _{1}(g^{2})+\chi
_{1}(g)^{2}-\sum_{\alpha ,\beta \in \Delta }\alpha .\beta e^{\alpha
+\beta }) \big{]}
\end{eqnarray*}%
where $\Delta $ denotes the roots of $\frak{g}$. $\chi _{2}(g)$ is an irreducible character for each of the Deligne exceptional Lie algebras. 
\end{theorem}
We expect that these methods can be extended to an analysis of the irreducible decomposition of $\dim V_n$ for $n\ge 3$ hopefully leading to a new understanding of Deligne's observations. 

\section{Griess Algebras}
\subsection{Genus Zero Constraints from Quadratic Casimirs}
Consider a VOA $V$ for which $d_1=\dim V_1=0$ and define a multiplication for $a,b \in V_2$ by $a\bullet b\equiv o(a)b=a_1b$. $V_2$ with multiplication $\bullet$ defines a commutative non-associative algebra, called a Griess algebra, for which the LiZ metric is a symmetric invariant bilinear form e.g. \cite{MN}. Let $V_2=\mathbb{C}\omega \oplus P_2$ where $P_2$ denotes
the space of weight $2$ primary vectors of dimension $p_2=\dim P_2$. Assume $p_2>0$ and let
$\{u_{\alpha}\}$ and $\{u^{\alpha }\}$ be a $P_2$ basis and LiZ dual basis and again define quadratic Casimir vectors by \cite{Mat,T1}
\begin{equation*}
\lambda ^{(n)}=u_{3-n}^{\alpha }u_{\alpha }\in V_n.
\end{equation*}
As before, we define a correlation function for $a,b\in P_2$
\begin{equation}
F(a,b;x,y)=\langle a,Y(u^{\alpha },x)Y(u_{\alpha },y)b\rangle.\label{eq:FabGriess}
\end{equation}
We then find that $F(a,b;x,y)$ is a rational function 
\begin{equation*}
F(a,b;x,y)=\frac{G(a,b;x,y)}{x^{4}y^{4}(x-y)^{4}},
\end{equation*}
where $G(a,b;x,y)$ is bilinear in $a,b$ and is a homogeneous, symmetric
polynomial in $x,y$ of degree $8$. We can next analyze the constraints on $F(a,b;x,y)$ that follow from $\lambda^{(4)}$ or $\lambda^{(6)}$ being Virasoro vacuum descendants to find results analogous to Theorems \ref{theorem:Killing} and \ref{theorem:d(C)} \cite{T1}:
\begin{theorem}
If $\lambda^{(4)}\in V_{4}(C,0)$ for $C\ne 0,-\frac{22}{5}$ then $G(a,b;x,y)$ is given by 
\begin{align*}
& \langle a,b\rangle \lbrack d_{2}{x}^{4}{y}^{4}+\frac{8d_{2}}{C}{x}^{3}{y}%
^{3}{(x-y)}^{2}+\frac{4d_{2}( 44-C) }{C( 5C+22) }{x}%
^{2}{y}^{2}{(x-y)}^{4} \\
& +2({x}^{2}+{y}^{2}){(x-y)}^{6}-{(x-y)}^{8}].
\end{align*} 
\end{theorem}
This implies that the $V_2$ trace form $Tr_{V_2}(o(a\bullet b))$ is given by 
\begin{equation*}
Tr_{V_2}(o(a\bullet b))=\frac{8(p_2+1)}{C}\langle a,b\rangle,
\end{equation*}
which being invertible, eventually implies $V_2$ is a simple Griess algebra. Furthermore for $C\neq C_{2,3}=0, C_{2,5}=-22/5, C_{2,7}=-68/7, C_{3,4}=1/2$ where $\det M_6(C,0)=0$ we have
\begin{theorem}
If $\lambda^{(6)}\in V_{6}(C,0)$ then
\begin{equation}
p_2(C)=\frac{1}{2}
\frac{(5C+22)(2C-1)(7C+68)}{C^2-55C+748}.\label{eq:d(C)Griess}
\end{equation}
\end{theorem}
This formula originally appeared in \cite{Mat} subject to stronger assumptions.
There are 37 rational values for $C$ for which $p_2$ is a positive
integer as follows \cite{T1,Ho1}:
\begin{eqnarray*}
-{\frac {44}{5}},8,{\frac {52}{5}},16,{\frac {132}{7}},20,
{\frac {102}{5}},{\frac {748}{35}},{\frac {43}{2}},22,
{\frac {808}{35}},{\frac {47}{2}},24,{\frac {170}{7}},
{\frac {49}{2}},{\frac {172}{7}},{\frac {152}{5}},{\frac{61}{2}},
{\frac {154}{5}},\\
{\frac {220}{7}},{\frac {63}{2}},32,{\frac {164}{5}},
{\frac{236}{7}},34,{\frac {242}{7}},36,40,{\frac {204}{5}},
44,{\frac {109}{2}},{\frac {428}{7}},68,{\frac {484}{7}},
{\frac {187}{2}},132,1496.
\end{eqnarray*}
\subsection{Genus Zero Constraints from Quadratic Casimirs}
We consider the genus one partition function $%
Z(q)=Tr(q^{L_{0}-C/24})$ and define a genus one correlation function:
\begin{eqnarray*}
Z((u^{\alpha },x),(u_{\alpha },y),q)=\sum_{n\geq 0}Z(\lambda^{[n]},q)(x-y)^{n-4}.
\end{eqnarray*}
Apply the Zhu reduction formula to the LHS we find:
\begin{eqnarray*}
Z((u^{\alpha},x),(u_{\alpha},y),q)=&& 
Tr(o(u^{\alpha })o(u_{\alpha })q^{L_{0}-C/24})
+P_2(x-y,q)Z(u^{\alpha}\bullet u_{\alpha},q)\\&&+p_2P_4(x-y,q)Z(q).
\end{eqnarray*}
Alternatively, applying Zhu reduction to $Z(\lambda^{[6]},q)$ for vacuum Virasoro descendant $\lambda$ to find:
\begin{theorem}\label{theorem:ZdeqnGriess}
If $\lambda^{(6)}\in V_{6}(C,0)$ then $Z(q)$ is the unique solution with leading form $Z(q)=q^{-C/24}(1+0+O(q^2))$ to the differential equation
\begin{eqnarray}                                                       
&&\big{[}
\left(q\frac{\partial}{\partial q}\right)^3
+6E_2(q)\left(q\frac{\partial}{\partial q}\right)^2\\ \nonumber
&&+\left(6E_2(q)^2-\frac {15}{124}E_4(q)(7C^2 +80C+152)  \right) q\frac{\partial}{\partial q}\\ \nonumber
&&-\frac {35}{248}CE_6(q)(5C^2+66C+144)\big{]}Z(q)=0.\label{eq:GriessZ}
\end{eqnarray}
$Z(q)$ is convergent for $0<|q|<1$.
\end{theorem}
Following the same line of argument as for Theorem \ref{theorem:ZWZW} we also have:
\begin{theorem}\label{theorem:GriessGen}
$V$ is generated by $V_2$.
\end{theorem}
We may recursively solve for $Z(q)$ to obtain $\dim V_n$ as a rational function in $C$. Enforcing that $\dim V_n$ be positive integral for $n\le 400$ restricts the 37 possible rational values of $C$ to the following 9 values:
\begin{description}
\item[$C=-\frac{44}{5},d_2=1$]  This can be realized as the VOA $L(-\frac{44}{5},0)\oplus L(-\frac{44}{5},2)$ with automorphism group $\mathbb{Z}_2$ formed by the Virasoro minimal model for $C_{3,10}=-\frac{44}{5}$ together with an irreducible module of highest weight two.\footnote{My thanks are due to Geoff Mason for pointing this out}.  
\item[$C=8,d_2=155$] This can be realized as the fixed point free lattice VOA $V_{L}^{+}$ (fixed under the automorphism lifted from the reflection isometry of the lattice $L$) for the rank 8 even lattice $L=\sqrt{2}E_8$. The automorphism group is $O^{+}_{10}(2).2$ \cite{G}. 
\item[$C=16,d_2=2295$] The VOA $V_{L}^{+}$ for the rank 16 Barnes-Wall even lattice $L=\Lambda_{16}$ whose automorphism group is
$2^{16}.O^{+}_{10}(2)$ \cite{S}.
\item[$C=23\frac{1}{2},d_2=96255$] This can be realized as the integrally graded subVOA of H\"ohn's Baby Monster Super VOA $VB^{\natural}$ whose automorphism group is the Baby Monster group $\mathbb{B}$ \cite{Ho2}.
\item[$C=24,d_2=196883$] This can be realized as the Moonshine Module $V^{\natural}$  constructed as a $\mathbb{Z}_2$ orbifolding of the Leech lattice VOA and whose automorphism group is the Monster group $\mathbb{M}$ and with $Z(q)=J(q)$ \cite{FLM}.
\item[$C=32,d_2=139503$] The partition function is that of a $\mathbb{Z}_2$  orbifolding of $V_{L}$ for an extremal self-dual lattice $L$ of rank 32 \cite{Ho2}. These lattices are  not fully classified. Amongst the known examples, it is not known if any satisfy the condition $\lambda^{(6)}\in V_{6}(C,0)$. 
\item[$C=32\frac{4}{5},d_2=90117$] No known construction.
\item[$C=33\frac{5}{7},d_2=63365$] No known construction.
\item[$C=40,d_2=20619$] The partition function is that of a $\mathbb{Z}_2$  orbifolding of $V_{L}$ for an extremal self-dual lattice $L$ of rank 40 \cite{Ho2}. These lattices are also not fully classified. Again it is not known if any examples satisfying the condition $\lambda^{(6)}\in V_{6}(C,0)$ exist.
\end{description}
Clearly, the last four examples deserve further investigation.

\subsection{Possible Deligne-like Tensor Product Decompositions?}
In the light of the discussion of Section 2, it is natural to examine the possibility of universal irreducible decompositions of tensor products of the Griess algebra. Consider the Virasoro decomposition of the weight three space
\begin{equation}
V_3=\mathbb{C}[L_{-3}\mathbf{1}]\oplus L_{-1}V_2\oplus P_3,
\label{eq:V3Decomp}
\end{equation}
where $P_3$ denotes the space of weight three primary vectors of dimension $p_3=\dim P_3\ge 0$. Then recursively solving the differential equation (\ref{eq:GriessZ}) for $Z(q)$ we find $p_2=\dim(P_2)$ as in (\ref{eq:d(C)Griess}) and
\begin{eqnarray*}
p_3 =\frac{31C(5C+22)(2C-1)(7C+68)(5C+44)}
{6(C^2-86C+1864)(C^2-55C+748)}.
\end{eqnarray*}
From Theorem \ref{theorem:GriessGen} we know that
$V$ is generated by $V_2=\mathbb{C}\omega \oplus P_2$ so that for $P_2$ 
basis $\{u^{\alpha}\}$ we have:
\begin{equation*}
V_3=\mathbb{C}[L_{-3}\mathbf{1},u_{-2}^{\alpha}\mathbf{1},
u_{0}^{\alpha}u_{-1}^{\beta}\mathbf{1}].
\end{equation*}
Since the Griess algebra involves the symmetric product $\bullet$ then $P_2$ is generated by the antisymmetric part of $\mathbb{C}[u_{0}^{\alpha }u_{-1}^{\beta }\mathbf{1}]$. This suggests that a Deligne-like
formula may also hold with 
\begin{equation}
\mathrm{Anti}\left( P_2\otimes P_2\right) =X\oplus Y,\label{eq:Anti}
\end{equation}
with $\dim \mathrm{Anti}\left(P_2\otimes P_2\right) =\frac{1}{2}p_2(p_2-1)$, $\dim X=p_3$ and where  $\dim Y$ is given by
\begin{equation*}
\frac{(5C+22)(2C-1)(7C+68)(5C+44)(3C^2-134C+136)(14C^2-553C-2796)}
{24(C^{2}-55C+748)^2(C^{2}-86C+1864)}.
\end{equation*}
This can be partly verified in four cases. For $C=-\frac{44}{5}$ the VOA has only two primary vectors: the vacuum $\mathbf{1}$ and a unique primary of weight 2. Thus (\ref{eq:Anti}) is trivial with $\dim X=\dim Y=0$. For $C=8,23\frac{1}{2}$ and $24$ for the groups $\mathbb{M}$, $\mathbf{B}$ and $O_{10}^{+}(2).2$ respectively one can check that $p_3=\dim X$ and $\dim Y$ are indeed the dimensions for irreducible representations of dimension $\chi_i(1A)$, in Atlas notation \cite{Atlas}, as shown:
\begin{equation*}
\begin{tabular}{|c|c|c|c|}
\hline
$C$ & $p_2$ & $p_3$ & $\dim Y$ \\ \hline
$-\frac{44}{5}$ & $1$ & $ 0$ & $ 0$ \\ \hline
$8$ & $\chi_2(1A)=155$ & $\chi_5(1A)= 868$ & $\chi_{11}(1A)= 11067$ \\ \hline
$23\frac{1}{2}$ & $\chi_3(1A)=96255$ & $ \chi_6(1A)=9550635$ & $\chi_{14}(1A)= 4622913750$ \\ \hline
$24$ & $\chi_2(1A)=196883$ & $\chi_3(1A)= 21296876$ & $ \chi_6(1A)=19360062527$ \\ \hline
\end{tabular}
\end{equation*}

\section{Higher Weight Constructions}
In general we may consider VOAs with a LiZ metric for which the primary vectors $P_k$ are of lowest weight $k\ge 3$. Exactly as before, we
can construct Casimir vectors from a $P_k$ basis $\{u_{\alpha }\}$ and LiZ dual basis $\{u^{\alpha }\}$ for $\alpha=1\ldots  p_k=\dim P_k$
\begin{equation*}
\lambda^{(n)}=u_{2k-1-n}^{\alpha}u_{\alpha }\in V_n.
\end{equation*}%
Then we find that provided $\lambda^{(2k+2)}$ is a vacuum Virasoro descendant then
$Z(q)$ is uniquely determined by a differential equation of order $k+1$ and that $V$ is generated by $V_k$ \cite{T2}. 

We conclude with a number of examples. For $k=3$ we find
\begin{equation*}
p_3={\frac{(5C+22)(2C-1)(7C+68)(5C+3)(3C+46)}
{-5{C}^{4}+703{C}^{3}-32992{C}^{2}+517172C-3984},}
\end{equation*}
with zeros of the Kac determinant $\det M_8(C,0)$. For $C=48$ we find $p_3=42987520$ and, in general, the partition function is the same as the conjectured extremal VOA of H\"ohn \cite{Ho2} arising in Witten's recent speculations on a possible relationship between three dimensional black holes and extremal VOAs of central charge $24n$ for $n\ge 1$ \cite{Wi}. For $k=4$ we find $p_4$ is given by
\begin{equation*}
\frac{5}{2}\frac{(2C-1)(7C+68)(5C+3)(3C+46)(11C+232)(C+10)}
{(C-67)(5C^4-1006C^3+67966C^2-1542764C-12576)},
\end{equation*}
with zeros of $\det M_{10}(C,0)$. For C=72 we obtain $p_4=2593096792$, and, in general, the partition function for Witten's conjectured extremal VOA. Lastly, for $k=5$ we find $p_5=-\frac {q(C)}{r(C)}$ where
\begin{eqnarray*}
q(C)=&&(7C+68)(2C-1)(3C+46)(5C+3)(11C+232)\\ 
&&.(13C+350)(7C+25)(5C+126)(10C-7)\\
r(C)=&&-363772080000+25483483057200C-37323519053016C^2\\
&&-7407871790404C^3+484484459322C^4-11429170478C^5\\
&&+132180881C^6-760575C^7+1750C^8,
\end{eqnarray*}
with zeros of $\det M_{12}(C,0)$. Interestingly, $p_5$ is \textbf{not} integral for $C=96$ so that we do not obtain the partition function for Witten's conjectured extremal VOA in this case.

\bibliographystyle{amsalpha}

\begin{thebibliography}{Atlas}
\bibitem[Atlas]{Atlas}J.H. Conway, R.T. Curtis, S.P. Norton, R.A. Parker and R.A. Wilson, \textit{Atlas of finite groups}, Oxford, Clarendon, 1985.

\bibitem[B]{B} R. Borcherds, \textit{Vertex algebras, Kac-Moody algebras and the Monster}, Proc.Natl.Acad.Sci.U.S.A. \textbf{83} (1986), 3068--3071.

\bibitem[BPZ]{BPZ} A. Belavin, A. Polyakov, and A. Zamolodchikov, \textit{Infinite conformal symmetry in two-dimensional quantum field theory}, Nucl.Phys. \textbf{B241} (1984), 333--380.

\bibitem[C]{C} A.M. Cohen, \textit{Some indications that the exceptional groups form a series}, CWI.Quart. \textbf{9} (1996), 51--59.

\bibitem[CdeM]{CdeM} A.M. Cohen and R. de Man, \textit{Computational evidence for Deligne's conjecture regarding exceptional Lie groups}, C.R.Acad.Sci.ParisS\'{e}r. I Math. \textbf{322} (1996), 427--432.

\bibitem[D]{D} P. Deligne, \textit{La s\'{e}rie exceptionnelle de groupes de Lie (The exceptional series of Lie groups)}, C.R.Acad.Sci.ParisS\'{e}r. I Math. \textbf{322}  (1996), 321--326.

\bibitem[DdeM]{DdeM} P. Deligne and R. de Man, \textit{La s\'{e}rie exceptionnelle de groupes de Lie II (The exceptional series of Lie groups II)}, C.R.Acad.Sci.ParisS\'{e}r. I Math. \textbf{323} (1996), 577--582.

\bibitem[DM]{DM} C. Dong and G. Mason,  \textit{Integrality of $C_2$-cofinite vertex operator algebras}, Int.Math.Res.Notices  \textbf{2006} Article ID 80468 (2006), 1--15. 

\bibitem[FHL]{FHL} I. Frenkel, Y-Z. Huang and J. Lepowsky, \textit{On axiomatic approaches to vertex operator algebras and modules}, 
Mem.Amer.Math.Soc. \textbf{104} (1993) no. 494.

\bibitem[FLM]{FLM} I. Frenkel, J. Lepowsky and A. Meurman, \textit{Vertex operator algebras and the Monster}, New York, Academic Press, 1988.

\bibitem[G]{G} R.L.Griess,  \textit{The vertex operator algebra related to $E_{8}$
with automorphism group $O^{+}(10,2)$}, in The Monster and Lie
algebras, (Columbus, Ohio, 1996), Ohio State University
Math.Res.Inst.Public. \textbf{7}, Berlin, de Gruyter, 1998.

\bibitem[Ho1]{Ho1}  G. H\"{o}hn, \textit{Conformal designs based on vertex operator algebras}, Adv.Math. \textbf{217} (2008), 2301--2335.

\bibitem[Ho2]{Ho2} G. H\"{o}hn, \textit{Selbstduale Vertexoperatorsuperalgebren und das
Babymonster}, Ph.D. thesis, Bonn.Math.Sch. \textbf{286} (1996), 1--85.

\bibitem[Hu]{Hu} K. Hurley, \textit{The space of graded traces for holomorphic vertex operator algebras of small central charge}, math.QA/0606282.

\bibitem[Ka]{Ka}  V. Kac, \textit{Vertex operator algebras for beginners}, University Lecture Series, Vol. 10, Boston, AMS, 1998.

\bibitem[KR]{KR} V. Kac  and A.K. Raina, \textit{Bombay lectures on highest weight representations of infinite dimensional Lie algebras}, Singapore, World Scientific, 1987.

\bibitem[Li]{Li} H. Li, \textit{Symmetric invariant bilinear forms on vertex operator algebras}, J.Pure.Appl.Alg. \textbf{96} (1994), 279--297.

\bibitem[LL]{LL} J. Lepowsky and H. Li, \textit{Introduction to
vertex operator algebras and their representations}, Birkh\"{a}user,
Boston, 2004.

\bibitem[LM]{LM} J.M. Landsberg and L. Manivel,  \textit{Triality, exceptional Lie algebras and the Deligne dimension formulas}, Adv.Math. \textbf{171} (2002), 59--85.

\bibitem[Mas]{Mas} G. Mason, \textit{2-Dimensional vector-valued modular forms}, Ramanujan J., to appear.

\bibitem[Mat]{Mat} A. Matsuo, \textit{Norton's trace formula for the Griess algebra of a vertex operator algebra with large symmetry}, Commun.Math.Phys. \textbf{224} (2001), 565--591.

\bibitem[MMS]{MMS} H. Maruoka, A. Matsuo and H. Shimakura, \textit{Trace formulas for representations of simple Lie algebras via vertex operator algebras}, unpublished preprint, 2005.

\bibitem[MN]{MN} A. Matsuo and K. Nagatomo, \textit{Axioms for a vertex algebra and the locality of quantum fields}, Math.Soc.Jap.Mem. \textbf{4} (1999).

\bibitem[MT1]{MT1} G. Mason and M.P. Tuite,  \textit{The genus two partition function for free bosonic and lattice vertex operator algebras}, arXiv:0712.0628.

\bibitem[MT2]{MT2}  G. Mason and M.P. Tuite, \textit{Torus chiral n-point functions for free boson and lattice vertex operator algebras}, Comm. Math. Phys. \textbf{235} (2003), 47--68.

\bibitem[MTZ]{MTZ} G. Mason, M.P. Tuite and A. Zuevsky,  \textit{Torus chiral n-point functions for $\mathbb{R}$ graded vertex operator superalgebras and continuous fermion orbifolds}, Commun.Math.Phys. \textbf{283} (2008), 305--342.

\bibitem[P]{P} J. Polchinski, \textit{String Theory,} Volumes I and II,
Cambridge University Press, Cambridge, 1998.

\bibitem[S]{S} H. Shimakura, \textit{The automorphism group of the vertex operator algebra $V_L^+$ for an even lattice L without roots}, 
J. Alg. \textbf{280} (2004), 29--57. 

\bibitem[T1]{T1} M.P. Tuite, \textit{The Virasoro algebra and some exceptional Lie and finite groups}, SIGMA \textbf{3} (2007), 008.

\bibitem[T2]{T2} M.P. Tuite, To appear.

\bibitem[Wa]{Wa} W. Wang,  \textit{Rationality of Virasoro vertex operator algebras}, Int. Math. Res. Notices. \textbf{71} (1993), 197--211.

\bibitem[Wi]{Wi} E. Witten, \textit{Three-dimensional gravity revisited}, preprint, arXiv:0706.3359.

\bibitem[Z]{Z} Y. Zhu,  \textit{Modular invariance of characters of vertex operator algebras}, J. Amer.Math.Soc. \textbf{9} (1996), 237--302.
\end{thebibliography}

\end{document}